\documentstyle[12pt]{article}
\input amssym.def
\input amssym.tex
\begin{document}

\thispagestyle{empty}
\pagestyle{myheadings}

\markright{\centerline{TORRES}}

\def\bbN{{\rm I}\!{\rm N}}
\def\bbR{{\rm I}\!{\rm R}}
\def\bbZ{{\rm Z}\!\!{\rm Z}}
\def\bbP{{\rm I}\!{\rm P}}
\def\bbH{{\rm I}\!{\rm H}}
\def\autx{{\rm Aut}(X)}
\def\ordsigma0{{\rm ord}(\sigma_0)}
\def\vsigma0{v(\sigma_0)}
\def\vtau{v(\tau)}
\def\fixsigma0{{\rm Fix}(\sigma_0)}
\def\ordtau{{\rm ord}(\tau)}
\def\fixtau{{\rm Fix}(\tau)}
\def\fixtaud{{\rm Fix}(\tau^d)}
\def\xsigma0{X/\langle \sigma_0 \rangle}
\def\ksigma0{k(X/\langle \sigma_0 \rangle)}
\def\ktau{k(X/\langle \tau \rangle)}
\def\gtau{g_{\tau}}
\def\gsigma0{g_{\sigma_0}}
\def\pitau{\pi_{\tau}}
\def\xtau{X/\langle \tau \rangle}
\def\xsigma0{X/\langle \sigma_0 \rangle}
\def\pisigma0{\pi_{\sigma_0}}
\def\pitau{\pi_\tau}
\def\pigamma{\pi_{J_\gamma}}
\def\lsigma0{\langle \sigma_0 \rangle}
\def\ltau{\langle \tau \rangle}
\def\k0{k_{\sigma_0}}
\def\k1{k_{\tau}}
\def\lj{\langle J_\gamma \rangle}
\def\vtaud{v(\tau^d)}

\begin{center}
\null\vskip2.5truecm

{\Large \bf BOUNDING THE ORDER OF AUTOMORPHISMS OF CERTAIN CURVES}
\vskip1truecm
{Fernando Torres }\\
{\small International Centre for Theoretical Physics}\\
{\small Maths. Group, P.O. Box 586, 34100  Trieste-Italy}\\
{\small feto@ictp.trieste.it}
\end{center}
\vskip0.5truecm
\begin{abstract}
We study upper bounds on the order of automorphisms of non-singular curves
$X$ satisfying
at least one of the following hypothesis: 1) $X$ is an $m$-sheeted
covering of exactly one non-singular curve of genus $\gamma$, where $m$ is
prime; 2) the center of the group of automorphisms of $X$ is non-trivial.
\end{abstract}

\noindent {\bf Notation.} Throughout this paper, by a curve we mean a
non-singular, irreducible
and projective algebraic curve defined over an algebraically closed field
$k$ of characteristic $p$. Let
$X$ be a curve and $P \in X$;
\vspace{-1pt}
\begin{list}
\setlenght{\rightmargin 0cm}{\leftmargin 0cm}
\itemsep=0.5pt
\item[$\bullet$] $k(X)$ and $\autx$ will denote, respectively, the field
of rational functions and the group of automorphisms of $X$. The symbol
${\rm div}_\infty(f)$ will stand for the polar divisor of $f\in k(X)$.
\item[$\bullet$] For $\tau \in \autx$, $\ordtau$, $\fixtau$ and
$\vtau$ will denote, respectively, the order, the set of fixed points
and the number of fixed points of $\tau$. $\k1$ and
$\gtau$ will denote, respectively, the
field of rational functions and the genus of the quotient curve
$\xtau$. $\pitau$ will denote the natural morphism $X\!\to\! \xtau$.
\item[$\bullet$] $H(P)$ and $G(P)$ will denote, respectively, the
Weierstrass semigroup and the set of gaps at $P$.
\end{list}
\medskip

\noindent {\bf Introduction.} In this paper we study upper bounds
on the order of automorphisms of curves satisfying at least
one of the following hypothesis.
\vspace{-1.5pt}
$$
\begin{array}{rll}
H_1(m,\gamma) &\! : &\ X\ {\rm is\ an\ }m-{\rm sheeted\ covering\ of\
 exactly\ one\ curve\ }\tilde X\ {\rm of\ genus\ }\gamma,\\
&\! &\ {\rm where\ } m\ {\rm is\ } {\rm prime}.\\[4pt]
H_2   &\! : &\ {\rm The\ center\ of\ }\autx\ {\rm is\ non-trivial}.
\end{array}
$$

\vspace{-1.5pt}

Let $X$ be a curve of genus $g\ge 1$ and let $\tau$ be an
automorphism of $X$.
We assume $\vtau\ge 1$ if $g=1$. It is known that the
order of $\tau$ satisfies
\vspace{-1pt}
\begin{equation}\label{ordtau}
\ordtau \le \left\{
\begin{array}{ll}
            2g+1   & {\rm if\ } \ordtau\ {\rm is\ prime}\\
        2(2g+1)     & {\rm otherwise,}
\end{array}\right.
\end{equation}

\vspace{-1pt}

\noindent except for some exceptional cases occurring for wildly
ramified extensions
$k(X)\!\mid\! k_{\tau}$ (Wiman [Wi], Harvey [Har], Singh
[S, Thms. 3.3, 3.3'], Stichtenoth [St, \S4]).
\medskip

Suppose that $X$ satisfies $H_1(m,\gamma)$. The following discussion follows
from Accola's [A, \S4] (see also [A2, Chapters 4, 5]). Let
$G(X\!\!\mid\!\!\tilde X)$ be
the group of covering transformations of $X\!\to\! \tilde X$, and let
$\tau \in \autx \setminus G(X\!\!\mid\!\!\tilde X)$. Then
$\tau$ induces an automorphism $\tilde \tau \in {\rm
Aut}(\tilde X)$ whose order is the smallest $\tilde n \in
\bbN$ such that $\tau^{\tilde n} \in G(X\!\mid\!\tilde X)$. If $\gamma
\ge 2$ or $v(\tilde \tau)\ge 1$ if $\gamma=1$, from (1) we have
upper bounds for $\tilde n$ and hence
for $\ordtau$. For instance if $\ordtau$ is a prime different from
$\#G(X\!\mid\!\tilde X)$, then
\vspace{-1pt}
\begin{equation}
\ordtau \le 2\gamma +1.
\end{equation}

\vspace{-1pt}

\noindent (see \S2). We remark that one can also obtain information about
$\#\autx$ because
$\autx / G(X\!\mid\! \tilde X) $ is isomorphic to a subgroup $H$ of ${\rm
Aut}(\tilde X)$. For example Accola (loc. cit.) used this to give an
explicit construction of curves admitting of only the identity as an
automorphism. On the other hand if $k(X)\!\mid\! k(\tilde X)$ is a Galois
extension, then
\vspace{-1pt}
$$
\autx = m \# H.
$$

\vspace{-1pt}

\noindent Thus if $X$ has many automorphisms, then either $X$ does not
satisfy $H_1(m,\gamma)$ with $\gamma \ge1$ (e.g. Hermitian curves, see
[St]), or if $X$ does, then $\gamma =0$ (e.g. the Klenian curve, see [Hur];
the curve $y^2=x^p+x$, see [Ro]), or $H$ has many automorphisms (see [Mac]).

We also remark that $H_1(m,\gamma)$ is satisfied if $X$ is an $m$-sheeted
covering of a curve of genus $\gamma$ and $g> 2m\gamma +(m-1)^2$. The
hypothesis on $g$ implies the uniqueness property of $H_1(m,\gamma)$ by
means of one of Castelnuovo's genus bound (see 1.1). The existence of an
$m$-sheeted covering from $X$ to a curve of genus $\gamma$ can be
characterized by means of the existence of certain Weierstrass semigroups
as well as the existence of certain linear series on $X$ (see [T]).
\medskip

Now suppose that $X$ satisfies $H_2$. Fix $\sigma_0$ in the center of
$\autx$ with $m:= \ordsigma0$ being a prime. Let $\tau \in \autx \setminus
\lsigma0$. We bound $\ordtau$ by using the data $(m,\gsigma0)$. As the
main consequence of $H_2$ we can ``pushdown" the data $(\ordtau, \vtau)$
on $X$ to $({\rm ord}(\tilde\tau), v(\tilde\tau))$ on $\tilde
X:=X/\lsigma0$, where $\tilde\tau$ is the pushdown of $\tau$ to
$\tilde X$. Moreover, $X/\langle \sigma_0,\tau\rangle$ is isomorphic to
$\tilde X/\langle \tilde\tau\rangle$, and $\vtau$ satisfies an equation
of type
\vspace{-1pt}
$$
\vtau = mu+f,
$$

\vspace{-1pt}

\noindent where
$u \in \bbN$ and $f=\#\fixsigma0 \cap \fixtau$.  In particular, if
$m\nmid \ordtau$ and ${\rm Fix}(\tau^d)={\rm Fix}(\tau)$, for $d\mid
\ordtau$, $d<\ordtau$ we find
\vspace{-1pt}
$$
2\gsigma0 -2 +u+f = \ordtau
(2g_{\tau_1}-2+u+f),
$$

\vspace{-1pt}

\noindent where $\tau_1:= \sigma_0\circ \tau$.
If $X$ also fulfils $H(m,\gsigma0)$ the above relation
improves (2) (see \S3).
\medskip

Typical examples of curves satisfying both the hypothesis above are
the  2-sheeted coverings having genus large
enough. Assume that $X$ is a $2$-sheeted covering of a curve of genus
$\gamma$, and let $J_\gamma$ be an
involution on $X$ whose orbits are the fibers of the $2$-sheeted
covering. Then
 $J_\gamma$ is unique provided $g> 4\gamma +1$ (Farkas
[F,
Corollary 2], Accola [A, Lemma 5]). Also in this case
$J_\gamma$ belongs to the center of $\autx$ (Farkas, [F, Thm. 2]; Accola
[A1, Application 4]). Furthermore Farkas (loc.cit.) showed that
\vspace{-1pt}
$$
\vtau \le 4\gamma +4
$$

\vspace{-1pt}

\noindent for $\tau \in \autx \setminus \langle J_\gamma \rangle$.
For the case of hyperelliptic curves ($\gamma =0$) of genus $g > 1$ it is
well known that all
the possibilities for $\vtau$  in $\{0,1,2,3,4\}$ occur and the unique
restriction
on $\ordtau$ is the Riemann-Hurwitz formula for $k(X)\!\mid\!k_{\tau}$ (cf.
Hurwitz [Hur]). However, if $\gamma \ge 1$, $g > 4\gamma+1$, and if we
assume $\vtau\ge 1$ for $\gamma=1$
the situation for both $\vtau$ and $\ordtau$ is different as we can see from
(2) and the above relation involving $u$ and $f$.
For instance it
was announced by Yoshida [Yo] that if $\gamma =1$ and $g> 5$, then the
possibilities for $(\ordtau, \vtau)$ are
\vspace{-2pt}
$$
(3,3), (3,5), (3,4), (3,2), (5,2), (5,3), (7,3), (9,2), (12,1), (8,2),
(6,1), (6,2), (4,4),
$$

\vspace{-2pt}

\noindent provided ${\rm Fix}( J_1)\cap \fixtau \not= \emptyset $.
\medskip

The
prototypes of our results are the following rather simple examples. They
also illustrate the methods used here.
\medskip

{\bf Example 1.} Let $X$ be a  hyperelliptic
curve of
genus $g>1$ defined over $k$ with $p \not= 2$. Let $\tau \in
\autx$ such
that $\ordtau>2$ is  prime. Assume that $k(X)\mid k_{\tau}$ is tamely
ramified. Set $f:=\#{\rm Fix }(J_0)\cap\fixtau$. Then $(\vtau,f) \in
\{(4,0), (3,1), (2,2)\}$.
\smallskip

 The hypothesis on $g$
implies that $J_0$ and $\tau$ commute  with each other. Hence
if $P \in \fixtau \setminus {\rm Fix}(J_0)$, then
$J_0(P) \in \fixtau$. Thus there exists $u\in \bbN$ such that
$\vtau =
2u +f$. Let $\tilde \tau$ be the
pushdown of $\tau$ to $X/\langle J\rangle$. We have that
$\ordtau = {\rm ord}(\tilde \tau)$ because $\ordtau$ is odd.
Hence the Riemann-Hurwitz formula applied to $\pi_{\tilde \tau}$ gives
\vspace{-2pt}
$$
-2+u+f = \ordtau(-2+u+f),
$$

\vspace{-2pt}

\noindent and so $u+f=2$, which establishes the example.
\medskip

{\bf Example 2.} Let $X$ be a $2$-sheeted
covering of an elliptic curve $\tilde X$. Suppose that the genus of $X$
satisfies $g>5$, and let $\tau$ and $k$ be as in Example 1.
Then $(\ordtau,\vtau) \not\in
\{(5,2),(3,2),(5,3),(7,3),(9,2)\}$.
\smallskip

Suppose that such a $\tau$ exists.
Let $\tilde \tau$ be the pushdown of $\tau$ to
$\tilde X$. Then $\ordtau = {\rm ord}(\tilde \tau)$ because $\ordtau$ is
odd. Now since $3$ is the only possible odd order for a non-trivial
automorphism of $\tilde X$ fixing a point, we reduce the example to analyze
the case $(\ordtau,\vtau)=(3,2)$.
With the notation from the above example we
have $\vtau = 2u+f=2$ and so $u+f\in \{1,2\}$. Consequently applying
Riemann-Hurwitz to $\pi_{\tilde\tau}$ we find
\vspace{-2pt}
$$
u+f=3(2\tilde g -2 + u+f),
$$

\vspace{-2pt}

\noindent where $\tilde g$ stands for the genus of $k(\tilde X)/\langle
\tilde \tau\rangle$. This is a contradiction.
\medskip

In particular, we see that not all the cases listed by
Yoshida can occur. In other words, the way as we bound $\ordtau$ gives
better results than his. We will also see that most of the known
results on automorphisms of hyperelliptic curves (e.g. Farkas-Kra [F-K;
III.7.11, V.2.13]) will emerge as simple corollaries of ours (see 5.1).
Yoshida and Farkas - Kra use Lewittes' results concerning
representations of the group of automorphisms as linear maps of
differential spaces (see [L]). To compute diagonal matrices, here one uses
the sequence of Weierstrass gaps at fixed points. Then, by means of the
character of the representation, one produces an equation $(*)$ involving
the
genus of the curve, the order of the automorphism and the number of its
fixed points. This equation and the Riemann-Hurwitz formula
imply restrictions for the order and the number of fixed points. When the
curve satisfies $H_2$ we obtain an analogous of $(*)$ by pushing down the
automorphism to an appropriated curve.
The advantage of this equation is that it does not
involve gaps sequence at fixed points.
\medskip

The contents of the paper are as follows. In \S1 we summarize the
results needed for the results stated here. We mainly based our
computations on one of Castelnuovo's genus bound (1.1), the
Riemann-Hurwitz formula (1.3) and on some results involving Weierstrass
semigroups (1.2).

In \S2 and \S3 we bound the order on automorphisms of
curves
satisfying hypothesis $H_1(m,\gamma)$ and $H_2$ respectively. In \S4
we consider necessary and sufficient conditions for automorphisms having
large number of fixed points. In 4.3 we improve Farkas' [F, Thm. 1].

In \S5 we specialize \S2 and \S3 to the case of double
coverings of curves. In 5.4 we consider
automorphisms of elliptic-hyperelliptic
curves. In 5.5 we deal with automorphisms of certain double coverings of
hyperelliptic curves, and we finish with 5.6 where we indicate how to
obtain results similar to those of 5.4 and 5.5 for certain double
covering of trigonal curves.

\bigskip

\noindent {\large \bf 1. Preliminary results}
\bigskip

\noindent {\bf 1.1 Castelnuovo's lemma ([C], [St1]).} Let $X$ be a curve
of genus $g$. Let $k_1$ and $k_2$ be two subfields of $k(X)$ with
 compositum equal to $k(X)$. Let $n_i$ (resp. $g_i$) be equal to the
degree (resp. genus) of $k_i$. Then
\vspace{-1pt}
$$
g \leqq n_1g_1 + n_2g_2 + (n_1 - 1)(n_2 -1).\quad \Box
$$

\vspace{-1pt}

\noindent {\bf 1.2 Remarks on Weierstrass semigroups}. Let $X$ be a
curve and $\tau \in \autx$ with $p\nmid \ordtau$.
\smallskip

{\bf (i)} If $P \in \fixtau$ and $h
\in \bbN$, then
\vspace{-1pt}
$$
\ordtau h \in H(P) \Leftrightarrow h \in H(\pitau(P)).
$$

\vspace{-1pt}

This is included in an implicit way in Kato's [K, p. 393]
(see also [T, Lemma 3.4]). Consequently ([Sch])
\vspace{-1pt}
$$
\gtau = \#\{\ell \in G(P): \ell \equiv 0\ ({\rm mod\ }\ordtau)\}.
$$

\vspace{-1pt}

{\bf (ii)} The above remark implies the following. Let $\sigma\in
\autx$ such that ${\rm ord}(\sigma) = \ordtau $ and ${\rm
Fix}(\sigma)\cap \fixtau \not=\emptyset $. Then
\vspace{-1pt}
$$
g_{\sigma} = \gtau.
$$

\vspace{-1pt}

{\bf (iii)} Let $P\in \fixtau$ and $\ell \in G(P)$ such
that $\ell \equiv 0$ (mod $\ordtau$). Since
$H(\pitau (P))\supseteq \{2\gtau,2\gtau+1,\ldots \}$, then (i) also implies
\vspace{-1pt}
$$
\ell \le (2\gtau-1)\ordtau.
$$

\vspace{-1pt}

In particular, if $\ordtau=2$ then $H(P)$ has $\gtau$ odd non-gaps
$\le 2g-1$. Moreover, let $U_1<\ldots<U_{\gtau}$ be such a non-gaps. Then
$U_1\ge 2g-4\gtau+1$ and
\vspace{-1pt}
$$
H(P)=\langle 2m_1,\ldots,2m_{\gtau},4\gtau+2,U_1,\ldots,U_{\gtau} \rangle,
$$

\vspace{-1pt}

\noindent where the $m_i$ are the first $\gtau$ positive non-gaps at
$\pitau(P)$ ([T, Lemmas 2.1, 2.3]).
\medskip

\noindent {\bf 1.3 The Riemann-Hurwitz formula.} Let $X$ be a
curve of genus
$g$, and $\tau \in \autx$. Assume $p\nmid
n:= \ordtau$. We will use the following
version of the Riemann-Hurwitz formula for $\pitau$ ([F-K, p. 274])
\vspace{-1pt}
$$
2g-2=n(2\gtau-2) + \mathop{\sum}\limits^{}_{d\mid
n,d<n}\varphi(n/d)v(\tau^d),
$$

\vspace{-1pt}

\noindent where $\varphi$ is the Euler function. In the formulae of \S3 we
will also use the number
\vspace{-1pt}
$$
\Lambda_\tau:= \mathop{\sum}\limits^{}_{d\mid n, d<
n}\varphi(n/d)(v(\tau^d)-\vtau).
$$

\vspace{-1pt}

\noindent The following
definitions allow us to have
a way of computing $\Lambda_\tau$ ([F-K, p.261]). The ramification set
of $\pitau $ can be partitioned into a disjoint
union of subsets $B_d$ with $d\mid n$ and $d<n$
where
\vspace{-1pt}
$$
B_1 = \fixtau ,\ \ \ B_d = \{P \in X: \tau^d(P)=P,\ \tau^i(P) \not= P\ {\rm
for}\ 0<i<d \}\ \ \ {\rm for}\ d>1.
$$

\vspace{-1pt}

\noindent Let $x_d=x_d(\tau):=\#\pitau(B_d)$. Then
\vspace{-1pt}
$$
\Lambda_\tau = \mathop{\sum}\limits^{}_{d\vert n,
1<d<n}(n-d)x_d.
$$
\noindent {\large \bf 2. Consequences of \boldmath $H_1(m,\gamma)$}
\bigskip

\noindent Let $X$ be a curve of genus $g$ satisfying $H_1(m,\gamma)$. Let
$\pi: X\!\to\! \tilde X$ be the $m$-sheeted covering of $X$ over a curve
of genus $\gamma$,
and let
$G=G(X\!\mid\!\tilde X)$ be the group of cover transformations of $\pi$. We
have $\#G =1$ or $\#G=m$ and the last case occurs if and only
if $k(X)\!\mid\!k(\tilde X)$ is a Galois extension. Let $\tau \in
\autx$. By the uniqueness of $\pi$
the pushdown $\tilde \tau$ of $\tau$ to $\tilde X$ is an automorphism of
$\tilde X$.
By means of the data $(m,\gamma)$ and by using (1), we
will set up upper bounds
on the order of $\tau \in \autx \setminus G$. If $\gamma=1$ we assume
$\vtau\ge 1$. (Since $\tilde\tau \circ \pi=\pi\circ\tau$, this implies
$v(\tilde\tau)\ge 1$.)

We have that ${\rm ord}(\tilde \tau)\mid \ordtau$, and ${\rm
ord}(\tilde \tau)$ is
the smallest positive integer $\tilde n$ such that $\tau^{\tilde n}\in G$.
Thus,
$$
{\ordtau \over {\rm ord}(\tilde \tau)}\mid \#G.
$$

\vspace{-3pt}

We
consider two cases.
\medskip

\noindent {\bf 2.1.} $ \ordtau = {\rm  ord}(\tilde
\tau)$.
(This is the case if $\#G=1$ or $m\nmid \ordtau$.) Here by (1) we have
$$
\ordtau \le \left\{
\begin{array}{ll}
            2\gamma+1   & {\rm if\ } \ordtau\ {\rm is\ prime}\\
        2(2\gamma+1)     & {\rm otherwise.}
\end{array}\right.
$$
\noindent {\bf 2.2.} $\ordtau \not= {\rm
ord}(\tilde \tau)$. Here we have $\ordtau =
{\rm ord}(\tilde \tau)m$, and hence $G = \langle \tau^{{\rm ord}(\tilde
\tau)}\rangle$. Thus (1) implies
\vspace{-1pt}
$$
\ordtau \le \left\{
\begin{array}{ll}
            (2\gamma+1)m   & {\rm if\ } {\ordtau\over m}\ {\rm is\ prime}\\
        2(2\gamma+1)m     & {\rm otherwise.}
\end{array}\right.
$$

\vspace{-1pt}

Once we know that $X$ admits an
$m$-sheeted covering
over a curve of genus $\gamma$, we have the following criterion for the
uniqueness of this covering:
\medskip

\noindent {\bf 2.3. Claim.} Let $X$ be a curve of genus $g$, $m$ a prime and
$\gamma$ a natural. If
\vspace{-1pt}
$$
g>2m\gamma + (m-1)^2,
$$

\vspace{-1pt}

\noindent then $X$ admits at most one $m$-sheeted covering over a curve
of genus $\gamma$.
\medskip

{\bf Proof.} Direct application of Castelnuovo's genus bound (1.1).\quad
$\Box$
\medskip

\noindent {\large \bf 3. Consequences of \boldmath $H_2$}
\bigskip

\noindent Let $X$ be a curve of genus $g$. Throughout this section we fix
$\sigma_0 \in \autx$ with $m:= \ordsigma0$ being a prime. Let $\tilde
X:= X/\lsigma0$ and $\tau
\in \autx$. Then the pushdown $\tilde\tau$ of $\tau$ to $\tilde X$
defines an automorphism of $\tilde X$ and thus
we can apply \S2 to $\pisigma0$. (We assume $\vtau\ge1$ if $\gsigma0=1$.)
However, if
\vspace{-2pt}
$$
\sigma_0\ {\rm  belongs\ to\ the\ center\ of\ } \autx\quad {\rm and}\quad
p\nmid \ordsigma0\ordtau,
$$

\vspace{-2pt}

\noindent  we can obtain more
precise
information on $\ordtau$.

The hypothesis on the center implies that $X/\langle
\sigma_0,\tau\rangle$ is isomorphic to $\tilde X/\langle \tilde
\tau\rangle$. Then by means of
the following particular equations for the number of fixed points, we can
pushdown the data $(\ordtau,\vtau)$ on $X$ to the data $({\rm
ord}(\tilde\tau),v(\tilde
\tau))$ on $\tilde X$.
\medskip

\noindent {\bf 3.1.}
Let $\tau \in \autx$ and set $n := \ordtau$. For $d\mid
n$ let
$$
f_d:= \# \fixsigma0 \cap \fixtaud.
$$
\noindent For $P \in
\fixtaud$,  $H_2$ implies
 $\{\sigma_0(P),\ldots,{ \sigma_0}^{m -1}(P)\} \subseteq \fixtaud$.
Since $m$ is prime we have $m$ points in the above set, unless $P \in
\fixsigma0$. Consequently, there exists a non-negative integer
$u_d=u_d(\sigma_0,\tau)$ such that
\vspace{-1pt}
\begin{equation}
v(\tau^d) = m u_d + f_d.
\end{equation}

\vspace{-1pt}

\noindent {\bf 3.2. Bounding the order I.}\quad Suppose
\vspace{-1.5pt}
$$
\lsigma0 \not\subseteq \ltau.
$$

\vspace{-1.5pt}

\noindent Here ${\rm ord}(\tilde\tau)=\ordtau$.
Thus if $\gsigma0\ge 2$ or $\vtau\ge 1$ for $\gsigma0=1$, 2.1
implies
$$
n=\ordtau \le \left\{
\begin{array}{ll}
2\gsigma0+1 & {\rm if\ } \ordtau\ {\rm is\ prime}\\
2(2\gsigma0+1) & {\rm otherwise}.
\end{array}
\right.
$$

Moreover, for any $\gsigma0$  we have the following
\medskip

\noindent {\bf 3.2.1. Lemma.}
\vspace{-1pt}
$$
\vtau \le {2m\gsigma0 \over n-1}+2m.
$$

{\bf Proof.} From 1.3 we have $n(2\gtau-2)+(n-1)\vtau \le 2g-2\ (*)$.
Now since $m$ is prime we have $k(X)=k_{\sigma_0}k_{\tau}$ and thus
the lemma follows from $(*)$ and Castelnuovo's genus bound (1.1).\quad $\Box$
\medskip

{\bf Remark.} The proof above only uses an inequality from the
Riemann-Hurwitz formula for $\pitau$.
Hence the lemma is also valid when $p\mid \ordsigma0 \ordtau$.
\medskip

Next we will improve the upper bound on $n$. We consider two cases
according as $m\nmid n$ or $m\mid n$.
\medskip

\noindent {\bf 3.2.2. \boldmath $m\nmid n$.}
Here we have $v(\tilde\tau)=u_d+f_d$ for $d\mid n$, $d<n$. Let
$\tau_1:=
\sigma_0\circ \tau$. Then
$\langle\sigma_0,\tau\rangle=\langle\tau_1\rangle$, and by  applying
the Riemann-Hurwitz formula to
$k_{\sigma_0}\!\mid\!k_{\sigma_0}/\langle \tilde\tau\rangle$
we find
\vspace{-1.5pt}
\begin{equation}
2\gsigma0 -2 + u_1 + f_1 = n(2g_{\tau_1}-2 + u_1+ f_1) +
\mathop{\sum}\limits^{}_{d\mid n, d<n}\varphi(n/d)(u_d+f_d- u_1-f_1).
\end{equation}

\vspace{-1.5pt}

\noindent In particular this implies
\vspace{-1.5pt}
\begin{equation}
\vtau = {2m\gsigma0 \over n-1} + 2m -{2mng_{\tau_1} + \Lambda_\tau +
(m-1)\mathop{\sum}\limits^{}_{d\mid n,d<n}\varphi(n/d)f_d \over n-1}.
\end{equation}

\vspace{-1.5pt}

\noindent ($\Lambda_\tau$ was defined in 1.3.) Equation (5) yields to
the following considerations.
\medskip

\noindent {\bf 3.2.2.1. \boldmath $\Lambda_\tau =0$.}\quad Here
(5) becomes
\vspace{-1.5pt}
\begin{equation}
2\gsigma0 -2 + u_1 + f_1 = n (2g_{\tau_1}-2 + u_1 + f_1).
\end{equation}

\vspace{-1.5pt}

\noindent {\bf (i)} Suppose
\vspace{-1.5pt}
$$
2\gsigma0
-2 + u_1 + f_1 =0.
$$

\vspace{-1.5pt}

\noindent Thus either $\gsigma0 =1$, $u_1+f_1 =0$, or
$\gsigma0 =0$ and $u+f_1 =2$. In the first case we have
$
g-1 = n(\gtau -1),
$
\noindent and in the second case
\vspace{-1.5pt}
\begin{list}
\setlenght{\rightmargin 0cm}{\leftmargin 0cm}
\itemsep=0.5pt
\item[(1)] $g-1+m = n(\gtau +m-1)$, $\vtau = 2m$, $f_1 =0$,
\item[(2)] $2g+m-1 = n(2\gtau +m-1)$, $\vtau = m+1$, $f_1 = 1$, or
\item[(3)] $g = n\gtau$, $\vtau = f_1 =2$.
\end{list}

\vspace{-1.5pt}

\noindent From (5) we notice that if $\gsigma0=0$, then
$u_1+f_1 =2 \Leftrightarrow \Lambda_\tau =0$. Consequently
the above three statements
generalizes a result on
automorphisms of prime order on hyperelliptic curves stated in [F-K,
V.2.13].
\medskip

\noindent {\bf (ii)} Now suppose
\vspace{-1.5pt}
$$
2\gsigma0 -2 + u_1 + f_1 \not= 0.
$$

\vspace{-1.5pt}

\noindent Here we have $\gsigma0 \ge 1$ and $\gsigma0=1 \Rightarrow
\vtau \ge 1$. Hence $n\le 2\gsigma0+1$ for $n$ prime. This
upper bound fulfils for any $n$:
\medskip

\noindent {\bf Claim. (1)}\quad $n\le 2\gsigma0 +1$.
\smallskip

\noindent {\bf (2)}\quad $n=2\gsigma0
+1$\ $\Leftrightarrow$\ $g_{\tau_1}=0$ and $u_1+f_1 =3$. Consequently
$$
(\vtau, f_1) \in
\{(3m,0),(2m+1,1),(m+2,2),(3,3)\}.
$$

{\bf Proof.} From (6) we have
$
n \le 2\gsigma0 -2 +u_1 + f_1$, and we can also write
\begin{equation}
u_1+f_1 = {2\gsigma0 -2ng_{\tau_1} \over n-1}+2.
\end{equation}
\noindent Then
$
n(n-1) \le 2n\gsigma0-2ng_{\tau_1}$, which implies statement (1).
To prove (2) we notice that $n=2\gsigma0+1$ implies $g_{\tau_1}=0$ except
for $\gsigma0=1$, $u_1+f_1=0$ (eliminated by our assumption). In fact if
$g_{\tau_1}\ge 1$, by (7) we have $u_1+f_1 \le
(n-3)/(n-1)$, and hence $u_1+f_1=0$. Thus once again by (8) we have $n=3$
and so $\gsigma0=1$.\quad$\Box$
\medskip

\noindent {\bf Remarks. (i)} Equation (7) implies
\vspace{-1pt}
$$
n< 2\gsigma0 +1 \quad \Rightarrow \quad n \le \gsigma0 +1,
$$

\vspace{-1pt}

\noindent unless the case
$g_{\tau_1}=1$, $u_1+f_1 =1$ where $n$ could be
$2\gsigma0 -1$.
\smallskip

{\bf (ii)}
\vspace{-1pt}
$$
n=\gsigma0 +1\quad \Leftrightarrow\quad u_1+f_1=4\ {\rm
and\ } g_{\tau_1}=0.
$$

\vspace{-1pt}

Consequently
\begin{list}
\setlenght{\rightmargin 0cm}{\leftmargin 0cm}
\itemsep=0.5pt
\item[(ii.1)] \quad $(\vtau,f_1)\in \{(4m,0), (3m+1,1), (2m+2,2),
(m+3,3), (4,4)\}$.
\item[(ii.2)] \quad $(\gsigma0 +1)m \in H(P)$ for $P\in
\fixsigma0\cap\fixtau$.
\end{list}

Let $n=\gsigma0+1$. Then (7) gives $u_1+f_1= 4 - 2ng_{\tau_1}/(n-1)\
(*)$. If
$g_{\tau_1}\ge 1$, then $u_1+f_1 \le 1$. If $u_1+f_1=1$ (resp. 0) $(*)$
gives
$\gsigma0=2$ (resp. 1). The first case leads to a contradiction because an
automorphism of prime order cannot have just one fixed point ([Gue],
[F-K, Thm. V.2.11]), and the second one is eliminated by hypothesis.
\bigskip

\noindent {\bf 3.2.2.2. \boldmath $\Lambda_\tau >0$.} By 1.3,
$\Lambda_{\tilde\tau}=\sum_{d\mid n, 1<d<n}^{}(n-d)x_d \ge n/2$. Hence
by (4) and 2.1, we obtain
\medskip

\noindent {\bf Claim. (1)}  Let $g_{\tau_1}\ge 1$. Then
$$
n \le \left\{
\begin{array}{ll}
   2(2\gsigma0 +1)/7      &    {\rm if\ } u_1+f_1 \ge3\\
       4 \gsigma0/5         &    {\rm if\ } u_1+f_1 =2\\
    2(2\gsigma0 -1)/3         &    {\rm if\ } u_1+f_1 =1\\
    4(\gsigma0-1)       &    {\rm if\ } u_1+f_1 =0.
\end{array}
\right.
$$
\noindent {\bf (2)} Let $g_{\tau_1} =0$ (recall that
$\gsigma0=1\Rightarrow u_d+f_d\ge 1$). Then
$$
n\le \left\{
\begin{array}{ll}
   4(\gsigma0 +1)/5           &    {\rm if\ } u_1+f_1 \ge 4\\
   2(2\gsigma0 +1)/3          &    {\rm if\ } u_1+f_1 =3\\
   4\gsigma0            &    {\rm if\ } u_1+f_1 =2\\
   2(2\gsigma0+1)                 &    {\rm if\ } \gsigma0
\ge 1\ {\rm and\ } u_1+f_1\le 1.
\end{array}
\right.
$$
\rightline{$\Box$}

\noindent {\bf 3.2.3. \boldmath $m\mid n$.} Unlike the above case here
$\langle \sigma_0, \tau\rangle$ is
not cyclic. However, we can use 3.2.1 and 2.1 to bound $n$. We find
\vspace{-1.5pt}
$$
n\le \left\{
\begin{array}{ll}
\gsigma0+1  & {\rm if\ } \vtau\ge 4 \\
2\gsigma0 + 1 &  {\rm if\ } \vtau=3 \\
2(2\gsigma0+1)  & {\rm if\ } \gsigma0  \ge 1,\ \vtau\le 2.
\end{array}\right.
$$

\vspace{-1.5pt}

\noindent {\bf 3.2.3.1. Remarks (i).} Suppose $g>2m\gsigma0+(m-1)^2$, let
$d\mid n$, $d<n$. If $f_d =
\#\fixsigma0\cap\fixtaud \ge 1$, then $m\nmid {n\over d}$.
\smallskip

Indeed, if $d\mid {n\over m}$, then $\tau^{n/m}= \tau^{dr}$ for
certain $r\in \bbN$. If $f_d \ge 1$, by 1.2 (ii) the genus of
$\tau^{n/m} = \gsigma0$ and hence by $H_1(m,\gsigma0)$ (the hypothesis on
$g$ implies this; see 2.3) we have $\lsigma0 \subseteq \ltau$, a
contradiction.
\smallskip

{\bf (ii)} By the Riemann-Hurwitz
formula applied to $k_{\tau^{n/m}}\!\mid \! k_{\tau}$ we find
\vspace{-1.5pt}
$$
n(2\gtau -2)+ (n-m)\vtau + \Lambda^{(1)}_\tau = m(2g_{\tau^{n/m}}-2),
$$

\vspace{-1.5pt}

\noindent where $\Lambda^{(1)}_\tau:= \sum_{d\mid
n,d<n,d\not=n/m}^{}\varphi(n/d)(v(\tau^d)-\vtau)$. In particular for
$n>m$ we obtain
\vspace{-1.5pt}
$$
\vtau = {2mg_{\tau^{n/m}}-2n\gtau - \Lambda^{(1)}_\tau \over n-m}+2.
$$

\vspace{-1.5pt}

{\bf (iii)} Let $n=m^xq$ with $m\nmid q$. From the Riemann-Hurwitz
formula for $k_{\tau^m}\!\mid\!k_{\tau}$ we get
\vspace{-1.5pt}
$$
mq(2\gtau-2) + (m-1)\mathop{\sum}\limits^{}_{d\mid
q}\varphi(q/d)v(\tau^d) = q(2g_{\tau^m}-2).\qquad (*)
$$

\vspace{-1.5pt}

\noindent Thus
\vspace{-1.5pt}
$$
\vtau = {2g_{\tau^m}-2m\gtau-\Lambda^{(2)}_\tau\over m-1}+2,
$$

\vspace{-1.5pt}

\noindent where $\Lambda^{(2)}_\tau:= (m-1)\sum_{d\mid
q}^{}\varphi(q/d)(v(\tau^d)-\vtau)$. In particular if $n=m^x$ we have
\vspace{-1.5pt}
$$
\vtau = {2g_{\tau^m} - 2mg_{\tau}\over  m-1}+2.
$$

\vspace{-1.5pt}

\noindent Now if we compute $v(\tau^q)$ by using the above formula
and replacing it in $(*)$, we get
\vspace{-1.5pt}
$$
(m-1)\mathop{\sum}\limits^{}_{d\mid q,d<q}\varphi(q/d)v(\tau^d) =
2(m-1)(q-1)+2qg_{\tau^m}-2g_{\tau^{mq}}+2mg_{\tau^q}-2mq\gtau.
$$

\vspace{-1pt}

{\bf (iv)} Finally we consider the extension
$k_{\tau^q}\!\mid\!k_{\tau}$. Let $n=m^xq$ with $m\nmid q$. We find
\vspace{-1.5pt}
$$
\mathop{\sum}\limits^{}_{d\mid n,q\nmid d}\varphi(n/d)v(\tau^d)=
2m^x(q-1)+ 2m^xg_{\tau^q}-2n\gtau.
$$

\vspace{-1.5pt}

\noindent Thus if $q>1$ we have
\vspace{-1.5pt}
$$
\vtau = {2g_{\tau^q} - 2q\gtau-\Lambda^{(3)}_\tau\over q-1}
+2,
$$

\vspace{-1.5pt}

\noindent where $\Lambda^{(3)}_\tau:= \sum_{d\mid
n, q\nmid d}\varphi(n/d)(v(\tau^d)-\vtau).$
\bigskip

\noindent {\bf 3.3. Bounding the order II.}\quad Now suppose
\vspace{-2pt}
$$
\lsigma0 \subseteq \ltau.
$$

\vspace{-2pt}

\noindent Here we have ${\rm
ord}(\tilde\tau)=n/m$, $g_{\tau^{n/m}}=\gsigma0$,
$\vtau=f_1$ and hence by the Riemann-Hurwitz formula for $k(\tilde
X)\!\mid\!k(\tilde X)/\langle\tilde\tau\rangle$ (or by
3.2.3.1 (ii)) we find
\vspace{-1.5pt}
\begin{equation}
n(2\gtau-2)+ (n-m)\vtau+\mathop{\sum}\limits^{}_{d\mid {n\over m}, d<{n\over
m}}\varphi(n/d)(v(\tau^d)-\vtau)=m(2\gsigma0-2).
\end{equation}

\vspace{-1.5pt}

\noindent In particular we have:
\vspace{-1.5pt}
$$
\vtau \le {2m\gsigma0\over n-m}+2.
$$

\vspace{-1.5pt}

Then by (8) and 2.2 (recall our assumption: $\gsigma0=1
\Rightarrow \vtau\ge 1$) we obtain
\medskip

\noindent {\bf 3.3.1. Claim. (1)}\quad If $\gtau \ge 2$, then
\vspace{-1.5pt}
$$
n \le {m(2\gsigma0 -2 + \vtau)\over \vtau +2}.
$$

\vspace{-1.5pt}

\noindent {\bf (2)}\quad If $\gtau =1$ and $\vtau \ge 1$, then
\vspace{-1.5pt}
$$
n \le {m(2\gsigma0 -2 +\vtau)\over \vtau}.
$$

\vspace{-1.5pt}

\noindent {\bf (3)}\quad If $\gtau =0$ and $\vtau \ge 3$, then
\vspace{-1.5pt}
$$
n \le {m(2\gsigma0 -2 + \vtau) \over \vtau -2}.
$$

\vspace{-1.5pt}

\noindent {\bf (4)}\quad Suppose that $\gsigma0\ge 1$. If $\gtau=1$ and
$\vtau=0$, or $\gtau=0$ and $\vtau \le
2$, then
\vspace{-1.5pt}
$$
n \le 2(2\gsigma0 +1)m.
$$
\noindent {\bf 3.3.2. Remarks.} Let $\tau \in \autx$ and set $n:=\ordtau$.
\smallskip

{\bf (i)} Let $g>2m\gsigma0+(m-1)^2$. From 3.2.3.1 (i) we have the
following
criterion for the hypothesis of this section. If $m\mid n$ and $f_1\ge
1$, then $\lsigma0\subseteq \ltau$.
\medskip

Now suppose $n=m^xq$, with $m\nmid q$. For any
$\gsigma0$ and $\vtau\in \{1,2\}$ we can
use 1.2 (i), (iii) and 3.2.3.1 to obtain more information on $n$:
\smallskip

{\bf (ii)} Let $\vtau =2$. Then either
$g_{\tau^q}\not=0$, or
all the powers of $\tau$, different from $\tau^{n/m}$, whose order is prime
also have two fixed points.
\smallskip

Suppose that $g_{\tau^q}=0$. Then by 3.2.3.1 (iv) we have
\vspace{-2pt}
$$
\mathop{\sum}\limits^{}_{d\mid n, q\nmid d}\varphi(n/d)\vtaud = 2m^x(q-1),
$$

\vspace{-2pt}

\noindent which proves the remark.
\smallskip

{\bf (iii)} Let $\fixtau = \{P\}$ (hence $q>1$ by [Gue] or
[F-K,  V.2.11]). Suppose that $x\ge 2$ and let
$\bar q$ be the smallest proper divisor of $q$, with $\bar q=1$ if $q$ is
prime. If
$
m^{x-1} \bar q > 2g_{\tau^{n/m}} +1,
$
then $m^{x-1} \in G(P)$.
\smallskip

By the formula for $\vtau$ in 3.2.3.1 (ii) and the hypothesis
on $\bar q$ we have $\vtaud\le 2$. Then by using the last equation
in 3.2.3.1 (iii) we get
\vspace{-1.5pt}
$$
(m-1)\varphi(q)\le -2qg_{\tau^m}+2g_{\tau^{mq}}-2mg_{\tau^q}+2mq\gtau.
$$

\vspace{-1.5pt}

\noindent Consequently $g_{\tau^{mq}}\not=0$ and hence by 1.2 (i), $m^{x-1}
= {\rm ord}(\tau^{mq}) \in G(P)$.
\medskip

\noindent {\large \bf 4. Additional remarks.}
\bigskip

\noindent Throughout this section $X$ is a curve of genus $g$ and
$\sigma_0\in \autx$ with $m:= \ordsigma0$ being a prime.
\smallskip

\noindent {\bf 4.1.} Suppose that $X$ satisfies $H(m,\gsigma0)$, and let
$\tau \in \autx$ with $\ordtau =m$. Assume that $p\nmid m$. Then either
$\tau \in \lsigma0$,
$f_1=\#\fixsigma0 \cap \fixtau =0$, or $\vtau =mu_1\le {2m\gsigma0\over
m-1}+2m$.
\smallskip

This is an immediate consequence of 1.2 (ii) and 3.2.1.
\medskip

\noindent {\bf 4.2.} Suppose that $\sigma_0$ belongs to the center of
$\autx$. Let $\tau \in \autx\setminus \lsigma0$ and set $n:=\ordtau$. Assume
that $p\nmid mn$. This remark is
concerned with the maximum value for $\vtau$ in the following cases
\medskip

\qquad I. $m\nmid n$,\qquad {\rm and}\qquad II. $\lsigma0 \subseteq
\ltau$.
\medskip

\noindent I. From (5) we have $\vtau \le {2m\gsigma0\over
n-1}+2m-(m-1)f_1$. From this equation we also have
\begin{list}
\setlenght{\rightmargin 0cm}{\leftmargin 0cm}
\itemsep=0.5pt
\item[(1)]
\vspace{-5pt}
$$
\vtau = {2m\gsigma0 \over n-1}+2m\quad \Leftrightarrow\quad \Lambda_\tau =0,
\ \ \fixsigma0 \cap \fixtau = \emptyset,\ \ g_{\tau_1}=0.
$$

\vspace{-3pt}

\item[(2)] If $f_1\ge 1$, then $\vtau \le {2m\gsigma0\over n-1}+m+1$ and
\vspace{-3pt}
$$
\vtau = {2m\gsigma0\over n-1}+m+1\quad \Leftrightarrow\quad
\Lambda_\tau=0, \ \ \ f_d=1\ {\rm for\ }d\mid n, d<n,\ \ g_{\tau_1}=0.
$$
\end{list}
\noindent II. By 3.3 we have $\vtau \le {2m\gsigma0\over n-m}+2$, and
\vspace{-1.5pt}
$$
\vtau = {2m\gsigma0 \over n-m}+2\ \Leftrightarrow\ \gtau=0,\ \vtau =
v(\tau^d)\ {\rm for}\ \
d\mid n,\ d<n,\ d\not= {n\over m}.
$$

\vspace{-1.5pt}

The equivalence follows from 3.2.3.1 (ii) (recall that
$\gsigma0= g_{\tau^{n/m}}$).
\medskip

\noindent {\bf 4.3. On Farkas' result.} This remark is concerned with
Farkas' [F, Thm. 1].
\medskip

(1) Let $\tau \in \autx$ such that $\vtau > m(2\gsigma0+1)$. Then $
\tau\in \lsigma0$.

In particular, $\fixsigma0\subseteq \fixtau$ and if
$\tau \not=1$, then $\ordtau = \ordsigma0$.
\medskip

(2) If $g>2m\gsigma0 +(m-1)^2$, then $\lsigma0 $ is normal in $\autx$. In
particular if $m=2$, $\sigma_0$ belongs to the center of $\autx$.
\medskip

{\bf Proof.} If $\tau \not\in \lsigma0$, then either by 3.2.1 we have $\vtau
\le
2m(\gsigma0+2)$ or by 3.3, $\vtau \le 2\gsigma0+2$. This proves (1).
\smallskip

Now we prove (2). By 2.3,  $k_{\sigma_0}$ is the only
subfield of $k(X)$ having index $m$ and genus $\gsigma0$. Let $\tau \in
\autx$.
Then since ${\tau}^{-1}\circ\sigma_0\circ\tau $ also has order $m$ and genus
$\gsigma0$ we must have ${\tau}^{-1}\circ\sigma_0\circ\tau \in
\lsigma0$ and we are done.\quad $\Box$
\medskip

\noindent {\bf Note.} Let $\tau \in \autx \setminus \lsigma0$. The
proof of (1) above shows that if
$\vtau=m(2\gsigma0+2)$,
then $\lsigma0 \not\subseteq \ltau$. If in adition $\gsigma0\ge 1$, then
$\ordtau =2$.
\medskip

\noindent {\large \bf 5. Double coverings}
\bigskip

\noindent In this section we specialize our results to the case $m=2$.
In what follows
$\pi:X\!\to\! \tilde X$ is a double covering of curves of genus $g$ and
$\gamma$ respectively.
We assume
\vspace{-2pt}
$$
g > 4\gamma +1.
$$

\vspace{-2pt}

\noindent Hence there exists a unique involution $J_\gamma$ belonging
to the center of $\autx$ and such that $\tilde X = X/\langle
J_\gamma\rangle$ (see 2.3 and 4.3 (2)). This involution will take the
place of $\sigma_0$ in \S3. For $\tau \in \autx$ we write $n=\ordtau$.
We recall that $\vtau=2u_1+f_1$ where $u_1 \in
\bbN$ and $
f_1 = \#{\rm Fix}(J_\gamma) \cap \fixtau$.

Moreover
\vspace{-1.5pt}
$$
\vtau\le \left\{
\begin{array}{ll}
{4\gamma\over \ordtau -1}+4 & {\rm if\ }\langle J_\gamma\rangle
\not\subseteq\ltau \\
{4\gamma\over \ordtau-2}+2 & {\rm otherwise}
\end{array}
\right.
$$

\vspace{-1.5pt}

\noindent (see 4.2).
\medskip

To begin with
we can reprove
well known results on automorphisms of hyperelliptic
curves ([Hur], [F-K, Thm. V.2.13]).
\medskip

\noindent {\bf 5.1 Proposition.} Let $X$ be a hyperelliptic curve of genus
$g>1$. Let $\tau \in \autx \setminus \langle J_0\rangle$, set $n:= \ordtau$
and suppose $p\nmid 2n$. \smallskip

(i) If $n$ is odd then $\vtaud = \vtau$ for $d\mid n$, and there are
three possibilities:
\vspace{-2pt}
\begin{list}
\setlenght{\rightmargin 0cm}{\leftmargin 0cm}
\itemsep=0.5pt
\item[(1)] $g+1=n(\gtau+1)$, $v(T) = 4$, $f_1=0$,
\item[(2)] $2g+1=n(2\gtau+1)$, $v(T) = 3$, $f_1 =1$, or
\item[(3)] $g=n\gtau$, $v(T) = 2$, $f_1=2$.
\end{list}

In cases (2) and (3), $X/\ltau$ is
hyperelliptic.
\medskip

(ii) If $n$ be even, then $f_1\le 2$ and we have
\begin{list}
\setlenght{\rightmargin 0cm}{\leftmargin 0cm}
\itemsep=0.5pt
\item[(1)] $f_1= \vtau =2\ \Rightarrow\ \vtaud =2$ for $d\mid n$,
$d<n$, $d\not=n/2$.
\item[(2)] $f_1=\vtau =1\ \Rightarrow\ n=2q$ with $q$ being an odd.
\item[(3)] $f_1=0\ \Rightarrow\ \vtau\in \{0,2,4\}$,

\end{list}

In cases (1) and (2), $\xtau$ has genus 0.
\medskip

{\bf Proof.} (i) $n$ odd. Equation (4) becomes
\vspace{-2pt}
$$
2(n-1)= (u_1+f_1)(n-1)+\Lambda_{\tilde \tau} \qquad (*),
$$

\vspace{-2pt}

\noindent where $\tilde \tau$ is the pushdown of $\tau$ to $\tilde X$.
Thus $u_1+f_1\le 2$.
\smallskip

{\bf Claim.}\qquad $u_1+f_1=2$.
\smallskip

Now (i) is a particular case of 3.2.2.1 (i). The statement on
hyperellipticity follows from 1.2 (i).
\smallskip

{\bf Proof of the claim.}
Suppose $u+f_1 =1$. Then $(*)$ and 1.3 implies $n-1=\sum_{d\mid
n,1<d<n}^{}(n-d)x_d$ for certain $x_d \in \bbN$. Since $n-d \ge 2n/3$
then $d=1$, a contradiction.

A similar argument also shows that $u_1+f_1=0$ is impossible.
\smallskip

(ii) $n$ even. If $\langle J_0\rangle\not\subseteq \ltau$, then $f_1=0$
and $\vtau=2u_1\le 4$ (see 4.3). Let $\langle J_0\rangle
\subset \ltau$. Then by 3.3 $\vtau\le 2$. Now
equation (8) becomes
\vspace{-2pt}
$$
\vtau = 2-{\Lambda^{(1)}_\tau\over n-2}.
$$

\vspace{-2pt}

\noindent This implies (1). Now let $\vtau=f_1=1$ and set $n=2^xq$
with $q$ odd. If $x\ge 2$ by 3.3.2 (iii) we would have
$2^{x-1}\in G(P)$, a contradiction.
\quad $\Box$
\medskip

{\bf Remark.} The examples in [Hur], [Ho] and [F-K] show that all the
cases of the proposition occur.
\medskip

{}From now on we assume $\gamma>0$.
\medskip

\noindent {\bf 5.2. Proposition} Let $\pi: X\to \tilde X$ be a 2-sheeted
covering of
curves of genus $g$ and $\gamma$ respectively. Suppose $\gamma>0$ and
$g>4\gamma+1$. Let $\tau \in
\autx\setminus \langle J_\gamma\rangle$ such that $\vtau \ge 1$ if
$\gamma=1$. Set $n:=\ordtau$ and assume $p\nmid 2n$.
\medskip

(i) Let $n$ be odd.
\vspace{-2pt}
\begin{list}
\setlenght{\rightmargin 0cm}{\leftmargin 0cm}
\itemsep=0.5pt
\item[(1)] If $\vtaud=\vtau$ for all $d\mid n$, then $n\le 2\gamma+1$.
\item[(2)] Set $\tau_1:=J_\gamma\circ\tau$. If $\vtaud\not=\vtau$ for some
$d\mid n$, then
$$
n\le \left\{
\begin{array}{ll}
3(2\gamma+1)/11 & {\rm if\ } u_1+f_1\ge 3\\
3\gamma/4   &  {\rm if\ } u_1+f_1=2\\
3(2\gamma-1)/5  & {\rm if\ } u_1+f_1=1\\
3(\gamma-1)  & {\rm if\ } u_1+f_1 =0.
\end{array}
\right.
$$
\vspace{-7pt}

provided $g_{\tau_1}\ge 1$, and
$$
n\le \left\{
\begin{array}{ll}
3(\gamma+1)/4  & {\rm if\ } u_1+f_1\ge 4\\
3(2\gamma+1)/5  & {\rm if\ } u_1+f_1=3 \\
3\gamma   &    {\rm if\ } u_1+f_1=2\\
2(2\gamma+1)-1  &  {\rm if\ } u_1+f_1 \le 1.
\end{array}
\right.
$$
\vspace{-10pt}

otherwise.
\end{list}

(ii) Let $n$ be even.
\begin{list}
\setlenght{\rightmargin 0cm}{\leftmargin 0cm}
\itemsep=0.5pt
\item[(1)] If $\langle J_\gamma\rangle \not\subseteq \ltau$, then $f_1=0$ and
$$
n\le \left\{
\begin{array}{ll}
\gamma+1 & {\rm if\ } \vtau\ge 8\\
2\gamma &  {\rm if\ } \vtau=6 \\
2(2\gamma+1) & {\rm if\ } \vtau\le 4.
\end{array}
\right.
$$

\vspace{-2pt}

\item[(2)] If $\langle J_\gamma\rangle \subseteq \ltau$, then
$$
n\le \left\{
\begin{array}{ll}
2(2\gamma-2+\vtau)/(\vtau+2)  &  {\rm if\ } \gtau \ge 2\\
2(2\gamma-2+\vtau)/\vtau  & {\rm if\ } \gtau=1\ {\rm and\ }\vtau \ge 1\\
2(2\gamma-2+\vtau)/(\vtau-2)  &  {\rm if\ } \gtau =0\ {\rm and\ }\vtau
\ge3 \\
4(2\gamma+1)  &  {\rm if\ } \gtau=1\ {\rm and\ } \vtau=0,\ {\rm or\
}\\
 &\ \ \ \gtau=0\ {\rm and\ }\vtau\le 2.
\end{array}
\right.
$$
\end{list}

{\bf Proof.} If $n$ is odd (1) follows from 3.2.2.1 while (2)
follows from (4) and 2.1 (notice that here $\Lambda_{\tilde \tau}
\ge 2n/3$, where $\tilde\tau$ is the pushdown of $\tau$ to $\tilde
X$). If $n$ is even the bounds follow from 3.2.3 and 3.3. \quad $\Box$
\medskip

\noindent {\bf Remark.} Let $\tau$ be as in 5.2 and suppose $f_1 \ge 1$.
If $n\ge 2\gamma$, then $\xtau$ is hyperelliptic. If $n$ is even and
$n\ge 4\gamma$, then $\xtau$ has genus zero.

This remark follows from 1.2 (i), (iii).
\medskip

\noindent {\bf 5.3} Let
$\tau \in
\autx$ whose order $n$ is odd and assume $\vtau\ge1$ if $\gamma=1$.
Assume further that $\Lambda_\tau =0$ and $p\nmid 2n$.
We state some remarks on $\vtau$ in the case where $n$ is large enough. By
\S3.2.2.1 this means $n=2\gamma+1$ or $n=\gamma+1$. Thus $g_{\tau_1}=0$
($\tau_1=J_\gamma\circ\tau$).
\medskip

{\bf (i) \boldmath $n=2\gamma+1$.} By the claim of 3.2.2.1 we have the
following table.
\bigskip

\centerline{\begin{tabular}{|c|c|c|c|c|}     \hline
Case  &  I  & II  &  III  &  IV  \\  \hline
$\vtau$ & 6  &  5  &  4  &  3  \\   \hline
$f_1$  & 0  & 1  &  2  &  3 \\ \hline
\end{tabular}}
\bigskip

We notice that in Case II from the Riemann-Hurwitz formula for $\pitau$
we have $2g-4\gamma+1=(2\gtau+1)(2\gamma+1)$ and hence for $P\in {\rm
Fix }(J_\tau)\cap \fixtau$ we find
$$
H(P)=\langle 2m_1,\ldots 2m_\gamma, 4\gamma+2, 2g-4\gamma+1 \rangle,
$$
\noindent where $m_1,\ldots,m_\gamma=2\gamma$ are the first $\gamma$
positive non-gaps at $\pi(P)$ (see 1.2 (i) (iii)).

In Case III we have $2g-2\gamma+1=(2\gamma+1)(2\gtau+1) \in H(P)$ for $P
\in {\rm Fix}(J_\gamma)\cap \fixtau$.
\medskip

{\bf (ii) \boldmath $n=\gamma+1$.} By Remark (ii) of 3.2.2.1 we have
$2(\gamma+1)\in H(P)$ for $P\in {\rm Fix}(J_\gamma)\cap\fixtau$, and the
following table.

\bigskip
\centerline{\begin{tabular}{|c|c|c|c|c|c|}     \hline
Case  &  I  & II  &  III  &  IV  & V \\  \hline
$\vtau$ & 8  &  7  &  6&  5 & 4  \\   \hline
$f_1$  & 0  & 1  &  2  &  3 & 4 \\ \hline
\end{tabular}}
\bigskip

Moreover in Case II, $H(P)$ is as in Case II above; in Case III,
$2g-3\gamma+1 \in H(P)$; in Case IV, $2g-2\gamma+1 \in H(P)$; and in Case V,
$2g-\gamma+1\in H(P)$.
\medskip

\noindent {\bf 5.4. Elliptic-hyperelliptic curves.} Let $X$ be a
$2$-sheeted covering of an elliptic curve $\tilde X$.
Let $\tau \in \autx\setminus \langle J_1\rangle$ with $\vtau\ge 1$, and
set $n:=\ordtau$.
\medskip

\noindent {\bf I. \boldmath $n$ odd.} By 5.2 (i) we have $n=3$ and hence
the possibilities for $(\vtau,f_1)$ are those of
the table in 5.3 (i). However, as we will see in the remark below, the
cases $(6,0)$ and $(4,2)$ are not possible. Hence the possibilities are
listed below.
\bigskip

\centerline{\begin{tabular}{|c|c|c|}     \hline
Case  &  1  & 2        \\  \hline
$\vtau$ &    5  &   3  \\   \hline
$f_1$  &  1  &   3 \\ \hline
\end{tabular}}
\bigskip

\noindent {\bf II. \boldmath $n$ even.} By 5.2 (ii) $n \le 12$. Moreover
$n\not=10$ because $X$ does not admit automorphisms of order $5$ fixing a
point.
\medskip

\noindent II.1 $\langle J_1 \rangle \not\subseteq \ltau$. Here $f_1 =0$,
$n\le 6$ and $2\le\vtau =2u_1\le 4+4/(n-1)$. If $n=2$ then $\vtau\in
\{4,8\}$. The case $n=6$ is not possible. All these
statements will be proved in the remark below.
\medskip

\noindent II.2 $\langle J_1 \rangle \subseteq \ltau$. Here $f_1=\vtau$,
$n\le 12$ and we have $\gtau =0$ and $1\le \vtau \le 4$.
Equation (8) becomes
\vspace{-2pt}
$$
(n-2)\vtau +\Lambda^{(1)}_\tau = 2n.
$$

\vspace{-2pt}

Then
\smallskip

(1) $n=4 \Leftrightarrow \vtau=4$;
\smallskip

(2) $n= 6 \Rightarrow \vtau\in\{1,3\}$. In fact, we have
$\vtau +
v(\tau^2)=6$. By the odd case $v(\tau^2)\in \{5,3\}$ and hence the result.
Conversely $\vtau=3$ gives $(n-6)+\Lambda^{(1)}_\tau =0$ and hence $n=6$;
\smallskip

(3) $n=8 \Leftrightarrow \vtau=2$. If $n=8$, $2\vtau+v(\tau^2)=8$. By
3.2.3.1 (iii), $v(\tau^2)= 4-4g_{\tau^2}$. Since this number is
positive then we must have $v(\tau^2)=4$ and thus the result. The
implication ``$\Leftarrow$" follows from the other cases.
\smallskip

(4) $n=12 \Rightarrow \vtau=1$. Here we find
$2\vtau+v(\tau^2)+v(\tau^3)+v(\tau^4)=12\ (*)$ and
$v(\tau^3)=4-4g_{\tau^3}$.
If this number is 0 then $v(\tau^2)+v(\tau^4)=12$. But since
$v(\tau^2)\le 4$, $v(\tau^4)\le 5$ this is a contradiction.
Hence $(*)$ becomes
$2\vtau+v(\tau^2)+v(\tau^4)=8$. By the case $n=6$ we have
$v(\tau^2)+v(\tau^4)=6$ and thus the result.
\smallskip

We now summarize this discussion in a table.
\bigskip

\centerline{\begin{tabular}{|c|c|c|c|c|c|c|} \hline
Case  &  3  & 4  & 5  &  6  & 7 & 8   \\ \hline
n  & 2 & $4$& 4 & 6 & 8 & 12      \\ \hline
$\vtau$ &$\{4,8\}$ &$\{2,4\}$  &  4& \{1,3\} & 2&  1 \\ \hline
$f_1$  & 0  & 0 &  4  & \{1,3\} & 2 & 1 \\ \hline
\end{tabular}}
\bigskip

\noindent {\bf Remarks. (i)} Let $\tau\in
\autx\setminus\langle J_1\rangle$ with $n=\ordtau=3$. We
show that $(\vtau,f_1)\not\in\{(6,0), (4,2)\}$. Let $P_0\in {\rm
Fix}(J_1)$. Since $2\in H(\pi(P_0))$ we have $4\in H(P)$. Let $x\in k(X)$
such that ${\rm div}_\infty(x)=4P_0$.
Then by Castelnuovo's genus bound (Lemma 1.1)
we have
\vspace{-1pt}
$$
k(x)\subseteq k(\tilde X).
$$

\vspace{-1pt}

\noindent Let $r$ (resp. $s$) be the number of
fixed
points of $J_1$ which are (resp. not) totally ramified for
$k(X)\!\mid\!k(x)$. Let $2t$ be the number of points $P\in X\setminus {\rm
Fix}(J_1)$ such that $\pi(P)$ is totally ramified for $k(\tilde
X)\!\mid\!k(x)$.
\medskip

\noindent {\bf Claim.}\qquad $s$ is even, $r+s=2g-2$ and $r+t=4$.
\medskip

{\bf Proof.} By [M-P, 7.5], $k(X)\!\mid\!k(x)$ is a
Galois extension and then $s$ is even and $r+s$ is the number of fixed
points of $J_1$ which is $2g-2$. Now by the Riemann-Hurwitz formula we have
\vspace{-2pt}
$$
2g-2 = 4(-2)+ 3r + s +2t,
$$

\vspace{-2pt}

\noindent from where it follows that $r+t=4$. \quad $\Box$
\medskip

Let $B_\pi:= \pi({\rm Fix}(J_1))$. Since $\tau$ commute with $J_1$ we
have $\tilde\tau(B_\pi)=B_\pi$. In particular, the claim above shows that
at least $f_1\ge 1$. This eliminates the case $(\vtau,f_1)=(6,0)$.
Now suppose $(\vtau, f_1)=(4,2)$. Thus $g\equiv 2$ (mod 3) ($*$).
Moreover, we
can assume that $P_0$ is fixed by $\tau$ or the image under $\pi$ of its two
fixed points not fixed by $J_1$ is a totally ramified for $k(\tilde
X)\!\mid\!k(x)$. In the first case there exists
$Q\in \fixtau$, such that $\pi(Q)$ is a point of
ramification 2 for $k(\tilde X)\!\mid\!k(x)$. Let
$x^{-1}(x(Q))=\{Q,Q_1\}$. Then
$\tilde\tau (\pi(Q_1))=\pi(Q_1)$ and, since ${\rm gcd}(2,3)=1$, we have
$\tau(Q_1)=Q_1$. This implies $f_1\ge 3$. In the second case $2g-6\equiv
0$ (mod 3). This is impossible under $(*)$ and we also eliminated
$(\vtau,f_1)=(4,2)$.
\medskip

{\bf (ii)} Here we illustrate Case 1 and Case 6: $n=6$,
$\vtau=f_1=1$.
\smallskip

Let $g$ be a multiple of $3$, $q(x)$ a separable polynomial over $k$
such that $q(0)\not=0$, ${\rm deg}(q(x))= (g-3)/3$, and let
$p(x):=4x^3-A$, $A\in k^*$. Suppose
that ${\rm gcd}(p(x),q(x)))=1$. Now consider the curve $X$ whose plane
model is given by
\vspace{-2pt}
$$
z^4=p(x)(q(x^3))^2.
$$

\vspace{-2pt}

\noindent Clearly $X$ is a double covering of the curve $y^2=P(x)$
which is an elliptic curve. By Riemann-Hurwitz $X$ has genus $g$ and the
unique point $P_\infty$ over $x=\infty$ is a fixed point of $J_1$.
Moreover, the four points over $x=0$ are not fixed by $J_1$.
Let $\epsilon$ be a $3$-root of unity. The automorphism $\tau$ given by
\vspace{-2pt}
$$
(x,z)\longmapsto (\epsilon x,z),
$$

\vspace{-2pt}

\noindent has five fixed points, namely $P_\infty$ and the points
over $x=0$. Considering $J_1\circ\tau$ we illustrate  the case with $n=6$.
\medskip

{\bf (iii)} Here we illustrate Case 2 and Case 6: $n=6$, $\vtau=f_1=3$.
\smallskip

Let $g\equiv 1$ (mod 3), $q(x)$ a separable polynomial over $k$ such that
$q(0)\not=0$, ${\rm deg}(q(x)=(g-4)/3$ and $p(x)=4x^3-A$, $A\in k^*$.
Suppose that ${\rm gcd}(p(x),q(x))=1$. Let $X$ be the curve given by
\vspace{-2pt}
$$
z^4= p(x)(q(x^3))^2 x^2.
$$

\vspace{-2pt}

\noindent Then $X$ is a double covering of $y^2=p(x)$; the unique point
$P_\infty$ over $x=\infty$ and the two points over $x=0$ are fixed by $J_1$.
Let $\epsilon$ be a $3$-root of unity and $\tau$ the automorphisms
defined by
\vspace{-2pt}
$$
(x,z)\longmapsto (\epsilon^2 x,\epsilon z).
$$

\vspace{-2pt}

\noindent Then the fixed points of $\tau$ are $P_\infty$ and the two
ones over $x=0$. By considering $J_1\circ\tau$ we illustrate the  case
where $n=6$.
\medskip

{\bf (iv)} Let $n$ even and suppose that $\langle
J_1\rangle\not\subseteq\ltau$. Thus
$n\in \{2,4,6\}$. We prove that $n=2\Rightarrow \vtau\in\{4,8\}$ and $n=6$
is not possible. Let $n=2$. Hence
$g=2\gtau-1+\vtau/2$. Since ${\rm ord}(\tilde\tau)=2$ from the proof of
the claim in (i), $g-3$ must be an even number. This eliminates
$\vtau\in\{2,6\}$. Now let $n=6$. Riemann-Hurwitz gives
$g-3=6\gtau-8+\vtau+v(\tau^2)+v(\tau^3)/2$. By the cases $n=2$ and
$n=3$, $g-3$ is odd. This eliminates $n=6$.
\smallskip

Next we show that our results are sharp. Let $g\in \bbZ^+$ such that
$g\equiv
3$ (mod 4). Let $p(x)=x^4-A$ irreducible in $k[x]$, $q(x)$ a separable
polynomial of degree $(g-3)/4$ such that ${\rm gcd}(p(x), q(x))=1$. Let
$X$ be the curve given by
\vspace{-2pt}
$$
z^4=p(x)(q(x^4))^2.
$$

\vspace{-2pt}

\noindent Then $X$ is a double covering of the curve $y^2=p(x)$ and the four
points over $x=\infty$ are not fixed by $J_1$. $X$ admits the
automorphisms
$
\tau:\ (x,z)\longmapsto (\epsilon x,z)$ and $\tau^2$, where $\epsilon$
is a $4$-root of unity.

If $q(0)\not=0$ then $\tau^2$ has eight fixed points and $\tau$ has two
fixed points. If $q(0)=0$, $v(\tau^2)=\vtau =4$.
\medskip

{\bf (v)} Finally we consider $\langle J_1\rangle \subseteq \ltau$. We can
assume
$n\in \{4, 8, 12\}$. Case 5 is illustrated by the automorphism $(x,z)\to
(x,\epsilon
z)$, $\epsilon$ a $4$-root of unity, defined on any of the above curves.
Now we illustrate Case 7. Let $g\in \bbZ^+$ such that $g\equiv 4$ (mod 4).
Let $p(x)$, $q(x)$ be as in (iv) except that ${\rm deg}(q(x))=(g-4)/4$
and $q(0)\not=0$. The curve
\vspace{-2pt}
$$
z^4=p(x)(xq(x^4))^2,
$$

\vspace{-2pt}

\noindent is a double covering of $y^2=p(x)$ and admits the automorphism
$(x,z)\mapsto (\epsilon^2 x,\epsilon z)$, where $\epsilon$ is a
$8$-root of unity. It fixs the two points over $x=\infty$.

Now we consider Case 8. The curve in (ii) admits of the automorphism
$(x,z)\mapsto (\epsilon^4 x,\epsilon^3 z)$, where $\epsilon$ is a
12-root of unity. Its unique fixed point is the only one over $x=\infty$.
\medskip

{\bf (vii)} Suppose that
$X$ admits of an automorphism as listed in Cases 1-8 and let $x$ be as in
(i).
Then, since $k(X)\!\mid\!k(x)$ is Galois [M-P, 7.5], $X$ is defined by a
model plane as in the examples above. To see this one uses the well
known group structure of automorphisms fixing a point on elliptic curves
(see e.g. Silverman [Sil, Thm. 10.1]) and Kato's [K1, \S6] or
Garcia's [G, Lemma 7].
\medskip

\noindent {\bf 5.5. Certain double coverings of hyperelliptic curves.} Here
we generalize the previous example. We consider curves $X$ such
that there exists $P_0\in X$ with $4\in H(P_0)$ and such that $X$ is a
double covering of a curve $\tilde X$ of genus $\gamma\ge 2$. Let $g$ be the
genus of $X$, $J_\gamma$ an involution such that
$X/\langle J_\gamma\rangle =\tilde X$, and let $x\in k(X)$ such that
${\rm div}_\infty(x)=4P_0$. Then by Castelnuovo's inequality (1.1) and
1.2 (i) we obtain
\medskip

\noindent {\bf 5.5.1. Claim.} If $g> 2\gamma+3$, then

(1) $P_0 \in {\rm Fix}(J_\gamma)$.

(2) $\gamma=\{\ell\in G(P_0): \ell\ {\rm even}\}$.\quad $\Box$
\medskip

In particular $\tilde X$ is hyperelliptic. In what follows we assume
$g>4\gamma+1$. Let $\pi:X\to \tilde X$, $r$ (resp.
$s$) the number of fixed points of $J_\gamma$ which are (resp. not)
totally ramified in $k(X)\!\mid\!k(x)$ and
let $2t$ be the number of points $P\in X\setminus {\rm Fix}(J_\gamma)$ such
that $\pi(P)$ is totally ramified for $k(\tilde X)\!\mid\!k(x)$. By
[M-P, 7.5], $k(X)\!\mid\!k(x)$ is
Galois (here it is enough to assume $g\ge 3\gamma$). Hence as in the
claim of 5.4 Remark (i), we obtain
\medskip

\noindent {\bf Claim.}\quad $s$ is even,
$r+s=2g-4\gamma+2$ and $r+t=2\gamma+2$.\quad $\Box$
\medskip

Next we only consider automorphisms $\tau$
such that $n=\ordtau$ is prime and equals to either $2\gamma+1$ or
$\gamma+1$.
\medskip

\noindent {\bf Case \boldmath $n=2\gamma+1$.} As in Remark (i) of 5.1
here we
also eliminate the cases I and III of 5.3 (i). We illustrate the remaining
cases.

Let $g\in \bbZ^+$ such that $g\equiv -2$ (mod ($2\gamma+1$)). Let
$q(x)$ be a separable polynomial over $k$ such that $q(0)\not=0$ and
${\rm deg}(q(x))=(g-3\gamma)/(2\gamma+1)$. Let $p(x):= x^{2\gamma+1}-A$
be an irreducible polynomial over $k[x]$ such that ${\rm gcd}(p(x),
q(x))=1$. Considerer the curve
\vspace{-2pt}
$$
z^4=p(x)(g(x^{2\gamma+1}))^2.
$$

\vspace{-2pt}

\noindent Then the unique point $P_0$ over $x=\infty$ satisfies $4\in
H(P)$; $X$ is a double covering over the hyperelliptic curve $y^2=p(x)$;
the automorphism
\vspace{-2pt}
$$
(x,z)\longmapsto (\epsilon x,z),
$$

\vspace{-2pt}

\noindent where $\epsilon$ is a $(2\gamma+1)$-root of unity, has five
fixed points namely $P_0$ and the four ones over $x=0$. This illustrate
case II of 5.3 (i).

Now take $g\in \bbZ^+ $ such that $g\equiv \gamma$ (mod $(2\gamma+1)$),
$q(x)$,
$p(x)$ as above but with ${\rm deg}(q(x))=(g-3\gamma-1)/(2\gamma+1)$. Then
the curve given by
\vspace{-2pt}
$$
z^4=p(x)(q(x^{2\gamma+1}))^2x^2,
$$

\vspace{-2pt}

\noindent is a double covering of $y^2=p(x)$, $4\in H(P_0)$ where $P_0$
is the unique point over $x=\infty$. The automorphism
\vspace{-2pt}
$$
(x,z)\longmapsto (\epsilon^2 x,\epsilon z),
$$

\vspace{-2pt}

\noindent where $\epsilon$ is a $(2\gamma+1)$-root of unity, has three
fixed points: $P_0$ and the two ones over $x=0$. This illustrate Case IV
of 5.3 (i).
\medskip

\noindent {\bf Case \boldmath $n=\gamma+1$.} As in the previous
examples here one can show that Cases I, II and IV of 5.3 (ii)
are not possible. The remaining cases are illustrated below.

Let $x$ be a transcendental element over $k$, $g\in \bbZ^+$ such that
$g\equiv -2$ (mod ($\gamma+1$)). Let $q(x)$, $p(x)$ be separable polynomials
over $k[x]$ such that ${\rm deg}(q(x))= (g-3\gamma-1)/(\gamma+1)$,
${\rm deg}(p(x))=2$, $q(0)\not=0$, $p(0)\not=0$ and ${\rm gcd}(p(x),
q(x))=1$. Consider the curve $X$ defined by
\vspace{-2pt}
$$
z^4= p(x^{\gamma+1})(q(x^{\gamma+1})x)^2.
$$

\vspace{-2pt}

\noindent Then $X$ is a double covering of $y^2=p(x^{\gamma+1})$; the
unique point $P_0$ over $x=a$, for $a\in k$ a root of $p(x)$, satisfies
$4\in H(P_0)$, and it has two or four points over $x=\infty$ according as
$g-2\gamma+1$ is odd or even. Let $\epsilon$ be a $(\gamma+1)$-root
of unity and consider the automorphism of $X$ given by
\vspace{-2pt}
$$
(x,z)\longmapsto (\epsilon^2 x,\epsilon z).
$$

\vspace{-2pt}

\noindent Its fixed points are those over $x=0$ or $x=\infty$. This
illustrate cases III and V of 5.3 (ii).
\medskip

\noindent {\bf Remark.} Let $X$ be a hyperelliptic curve admitting of an
automorphism $\tau$ of order $n$. Then Hurwitz [Hur] showed that $X$
and $\tau$ can be
defined by $y^2=f(x^n)$ and $\tau: (x,y)\mapsto (\epsilon x, \pm y)$,
or by $y^2=xf(x^n)$ and $\tau: (x,y)\mapsto (\epsilon x,
\epsilon^{1/2} y)$, where $\epsilon$ is a $n$-root of unity.

Let $x$ be as in the above examples. Using the fact that
$K(X)\!\mid\!k(x)$ is Galois, the mentioned Hurwitz's results and Komeda's
[Ko,
\S4], one can show that the curves of this section admitting of an
automorphism
satisfying Cases III and V of 5.3 (ii), can be defined by a model plane as
those of the examples stated above.
\medskip

{\bf 5.6. Certain double coverings of trigonal curves.} To finish this
paper, let consider a curve $X$ admitting of a point
$P_0$ such that $6\in H(P_0)$ and such that $X$ is a double covering of a
curve $\tilde X$ of genus $\gamma$. Let $J_\gamma$ be an involution such
that $X/\langle J_\gamma\rangle = \tilde X$ and let $x\in k(X)$ such that
${\rm div}_\infty(x)= 6P_0$. As in 5.5.1 we have that $P_0\in {\rm
Fix}(J_\gamma)$ and  $\gamma$ is the genus of $\tilde X$, provided
$g>2\gamma+5$. In particular $\tilde X$ is trigonal. We further assume
\vspace{-2pt}
$$
g\ge\ {\rm max}(4\gamma+2,
2\gamma+\rho),
$$

\vspace{-2pt}

\noindent where $\rho :=\{ \ell \in G(P_0): \ell
\equiv 0\ ({\rm mod\ }3)\}$. Under this condition it follows from [M-P,
Thm 7.1] that $k(X)\!\mid\!k(x)$ is a Galois extension.
Then, it is not
difficult to see that
 $H(P_0)=\langle 6, 2g-4\gamma+1+2i_2+4i_4+4i_5,
4\gamma+4 -2i_1 -2i_4, 2g-4\gamma+1+2i_1+2i_5, 4\gamma+4-2i_2 -2i_5,
2g-4\gamma+1 +4i_1+4i_2+2i_4  \rangle$, where the $i_j's$ are
non-negative integers satisfying $i_1+i_2+i_4+i_5 = \gamma+1$,
$i_1+i_3+i_5 =
2r+1$ and $i_1+2i_2+i_4+2i_5 \not \equiv 0$ (mod $3$). Moreover, $X$
admits a model plane of type
\vspace{-2pt}
$$
z^6 =
\mathop{\prod}\limits^{5}_{j=1}\mathop{\prod}\limits^{i_j}_{i=1}(x-a_{ij})^j,
$$

\vspace{-2pt}

\noindent where the $a_{ij}'s$
are pairwise different
elements of $k$.

{}From this fact one can prove
results similar to 5.4 and 5.5. We leave these to the reader.
\bigskip

\centerline{\bf Acknowledgments}
\bigskip

\noindent Thank you very much to the International Centre for
Theoretical Physics, Trieste, and its ``piccolo" but charming
Mathematics Section. Special thanks to Profs. V. Brinzanescu and A.
Verjovsky for helpful discussions.
\bigskip

\centerline{\bf References}

\begin{description}
\small
\itemsep=0.5pt
\item{[A]}
Accola, R.D.: Strongly branched coverings of closed Riemann surfaces,
Proc. Amer. Math. Soc. {\bf 26}, 315--322 (1970).

\item{[A1]}
Accola, R.D.: Riemann surfaces with automorphism groups admitting
partitions, Proc. Amer. Math. Soc. {\bf 21}, 477--482 (1969).

\item{[A2]}
Accola R.D.: Topics in the theory of Riemann surfaces, Lecture notes in
Math. 1595, Springer-Verlag (1994).

\item{[C]}
Castelnuovo, G.: Sulle serie albegriche di gruppi di punti appartenenti ad una
curvea algebrica, Rendiconti della Reale Accademia dei Lincei (5) {\bf 15},
337--344 (1906).

\item{[F]}
Farkas, H.M.: Remarks on automorphisms of compact Riemann surfaces, Ann.
of Math. Stud. {\bf 78}, 121--144 (1974).

\item{[F-K]}
Farkas, H.M.; Kra, I.: Riemann surfaces, Grad. Texts in Math.
{\bf 71}, Springer-Verlag (second edition) (1992).

\item{[G]}
Garcia, A.: Weights of Weierstrass points in double coverings of curves
of genus one or two, Manuscripta Math. {\bf 55}, 419--432 (1986).

\item{[Gue]}
Guerrero, I.: On Eichler trace formulas, Modular functions in analysis and
number theory, Ed. T. Metzger, Univ. Pittsburgh (1980).

\item{[Har]}
Harvey, W.J.: Cyclic groups of automorphisms of a compact Riemann
surface, Quart. J. Math. Oxford, Ser. (2), {\bf 17}
, 86--97 (1966).

\item{[Ho]}
Horiuchi, R.: Normal coverings of hyperelliptic Riemann surfaces, J.
Math. Kyoto Univ. {\bf 19}, 3, 497--523 (1979).

\item{[Hur]}
Hurwitz, A.: \"Uber algebraische Gebilde mit eindeutigen Transformationen
in sich, Math. Ann. {\bf 41}, 403--442 (1893). Reprinted in Mathematische
Werke I., Birkh\"auser; Besel Verlag, 391--430 (1962).

\item{[K]}
Kato, T.: On the order of a zero of the theta function, Kodai Math. Sem. Rep.
{\bf 28}, 390--407 (1977).

\item{[K1]}
Kato, T.: Non--hyperelliptic Weierstrass points of maximal weight, Math.
Ann. {\bf 239}, 141--147 (1979).

\item{[Ko]}
Komeda, J.: On Weierstrass points whose first non-gaps are four, J. Reine
Angew. Math. {\bf 341}, 68--86 (1983).

\item{[L]}
Lewittes, J.: Automorphisms of compact Riemann surfaces, Amer. J. Math.
{\bf 85}, 734--752 (1963).

\item{[Mac]}
Macbeath, A.M.: On a theorem of Hurwitz, Proc. Glasgow Math. Assoc. {\bf
5}, 90--96 (1961).

\item{[M-P]}
Morrison, I.; Pinkham, H.: Galois Weierstrass points and Hurwitz
characters, Ann. of Math., {\bf 124}, 591--625 (1986).

\item{[Ro]}
Roquette, P.: Abscha\"atzung der Automorphismenanzahl von
Funktionenk\"orpern bei Primzahlcharakteristik, Math. Z. {\bf 117},
157--163 (1970).

\item{[Sch]}
Schoeneberg, B.: \"Uber die Weierstrass Punkte in den K\"orpern der
elliptischen Modulfunktionen, Abh. Math. Sem. Univ. Hamburg {\bf 17},
104--111 (1951).

\item{[Sil]}
Silverman, J.H.: The arithmetic of elliptic curves, Grad. Texts in
Math. {\bf 106}, Springer-Verlag, (1986).

\item{[S]}
Singh, B.: On the group of automorphisms of a function field of genus at
least two, J. Pure. Appl. Alg., {\bf 4}, 205--229 (1974).

\item{[St]}
Stichtenoth, H.: \"Uber die Automorphismengruppe eines algebraischen
Funktionenk\"orpes von Primzahlcharakteristik; Teil II: Ein spezieller
Typ von Funktionenk\"orpern, Arch. Math. (Besel), {\bf 24}, 615--631 (1973).

\item{[St1]}
Stichtenoth, H.: Die ungleichung von Castelnuovo, J. Reine Angew. Math. {\bf
348}, 197--202 (1984).

\item{[T]}
Torres, F.: On certain $N$-sheeted coverings of curves and
numerical semigroups which cannot be realized as Weierstrass semigroups,
to appear in Comm. Algebra.

\item{[Wi]}
Wiman, A.: Ueber die hyperelliptischen Curven und diejenigen vom
Geschlechte $p=3$, welche eindeutigen Transformationen in sich zulassen,
Bihang Till Kongl. Svenska Vetenskaps-Akademiens Handlingar (Stockholm 1895
- 96)

\item{[Y]}
Yoshida K.: Elliptic--hyperelliptic Weierstrass points and automorphisms
of compact Riemann surfaces, Poster ICM 94 (1994).

\end{description}
\vskip1.5truecm

\end{document}